\documentclass[13pt]{article}
\usepackage{geometry}                
 \newtheorem{theorem}{Theorem}
\newtheorem{corollary}[theorem]{Corollary}
\newtheorem{definition}{Definition}
\newtheorem{proposition}{Proposition}
\newtheorem{lemma}{Lemma}
\newtheorem{remark}{Remark}
\usepackage{graphicx}
\usepackage{amssymb}
\usepackage{epstopdf}
\usepackage{amssymb}
\title{ A Representation of Multiplicative Arithmetic Functions by Symmetric Polynomials}
\author{Trueman MacHenry and Kieh Wong}
\date{October 10, 2007}                                           

\begin{document}
\maketitle
\begin{center}
\textbf{ABSTRACT}
\end{center}

We give a representation of the classical theory of multiplicative arithmetic functions (MF)in the ring of   
symmetric polynomials. The basis of the ring of symmetric polynomials that we use is the isobaric basis, a  basis especially sensitive to the combinatorics of partitions of the integers. The representing elements are recursive sequences of Schur polynomials evaluated at subrings of the complex numbers. The multiplicative arithmetic functions are units in the Dirichlet ring of arithmetic functions, and their properties can be described locally, that is, at each prime number $p$.  Our representation is, hence, a local representation. One such representing sequence is the sequence of generalized Fibonacci polynomials. In general the sequences consist of Schur-hook polynomials. This representation enables us to clarify and generalize classical results, e.g., the Busche-Ramanujan identity, as well as to give a richer structural description of the convolution group of multiplicative functions.  It is a consequence of the representation that the MF's can be defined in a natural way on the negative powers of the prime $p$.  
\vspace{0.5cm}

\textbf{CONTENTS}	

0. Introduction

1. Ring of Isobaric Polynomials

2. The Companion Matrix and Core Polynomials

3 The Ring of Arithmetic Functions.

4. Relation Between Arithmetic Functions and the WIP-Module.

5. Examples .

6. Specially Multiplicative Arithmetic Functions.

7. Structure of the Convolution Group of Multiplicative Arithmetic Functions Reconsidered.

8. Norms.

9. WIP-Module Revisited.

\vspace{0.5cm}

0.  \textbf{INTRODUCTION} 

\vspace{0.25cm} In this paper we give a representation of the classical theory of multiplicative arithmetic functions (MF) in the ring of symmetric polynomials. The Dirichlet ring of arithmetic functions $\mathcal{A}^*$ is well known to be a unique factorization domain.
\cite{CE}. Its ring theoretic properties have been investigated in, e.g., \cite{R},\cite{R2}, \cite{HS}, \cite{CG}, \cite{TM}, \cite{MT}. The multiplicative arithmetic functions are units in this ring and their properties can be described locally, that is, at each prime number $p$,\cite{PJM},\cite{RS2} \cite{RV}. It is this local behaviour which we take advantage of to construct a representation in terms of a certain class of symmetric polynomials.  It turns out that it is advantageous to use as a basis of the ring of symmetric polynomials, the \textit{isobaric} basis, which we describe in Section 1.  This basis is especially sensitive to the combinatorics of partitions of the integers, \cite{MT}, \cite{MT2}, \cite{TM2}.  Henceforth,  we refer to the symmetric polynomials in this basis as the ring of isobaric polynomials.   

The link between the theory of symmetric polynomials and the theory of multiplicative, arithmetic functions is that of linear recurrences, especially the ideas contained in \cite{MW};  also see  \cite{JR}, \cite{L}, \cite{HM}.   The isobaric ring contains a certain submodule, the submodule of \textit{weighted isobaric polynomials}(WIP), which is generated as a $\mathbb{Z}$-module by the Schur-hook polynomials. This module has the property that it can be partitioned into sequences which are linear recursions \cite{MT},  and, in fact, every linear recursion can be obtained from these sequences.  It contains the sequence of Generalized Fibonacci Polynomials (GFP), and the sequence of Generalized Lucas Polynomials  (GLP), \cite{TM2}.  It turns out that each of these sequences,  when the indeterminates are evaluated over a subring of the complex numbers, is the evaluation of a local sequence of multiplicative functions, i.e., a multiplicative function at a prime $p$. Moreover, every MF is represented locally by such sequences, in fact, already by the GFP-sequence.  This fact brings all of the machinery of the isobaric ring to bear with respect to the convolution group of multiplicative functions.  The importance of linear recursions in the theory of MF is recognized in \cite{LP} and in \cite{JR};  however,  the connection between multiplicative functions and symmetric functions and the power of the isobaric notation to simplify and reveal basic facts about the structure of MF is not made explicit in these papers.  

The same machinery of linear recurrences in the isobaric setting was used in \cite{MW} to study number fields.  A consequence of the results in that paper implies a certain strong connection between the structure of number fields, the algebraic structure of multiplicative functions,and periodicity in the theory of recurrsion.  The connection between  $\mathcal{A}^*$ and the symmetric polynomials was exploited in \cite{TM} and in \cite{MT}; in the first of these papers it was used to prove that the group of multiplicative functions generated by the completely multiplicative functions is free abelian; in the second paper a constructive procedure using isobaric polynomials was given for embedding this group into its divisible closure---also see \cite{CG}. 
 
In Section 1, we define the \textit{weighted isobaric polynomials} (WIP's) and give a formula for them independent of recursion.  

In Section 2,  the notion of the \textit{core} polynomial is introduced and the \textit{infinite companion matrix} and its properties are described, (see also \cite{MW}, \cite{L2}. This infinite matrix extends the WIP-sequences in the negative direction (that is, provides negatively-indexed functions as well as positively-indexed ones).  This addition of the negatively indexed sequences of isobaric polynomials carries over to the MF's. Also,  in this section generating functions are provided for the isobaric polynomials (and their MF counterparts).  

In Section 3,  we discuss the ring of arithmetic functions, and prove an important classification theorem concerning their structure.

In Section 4,  the main theorem asserting the relation between multiplicative arithmetic functions and the WIP-module is proved.  In this section it is shown that for each MF $\alpha$, not only do we have its local representation in terms of GFP's, but in addition,  each column of the infinite companion matrix also determines an MF, which we might call a "satellite" function of  $\alpha$. Each element in any of these sequences is a Schur-hook polynomial. (The negatively indexed one's provide an extension of the idea of Schur polynomials giving rise to the problem of their combinatorial interpretation.) All of these Schur polynomials,  the negatively-indexed one's as well as the positively-indexed one's,  can be conveniently computed  using Jacobi-Trudi formulae in their isobaric form.

In Section 5,  we look at the theory of \textit{specially multiplicative }arithmetic functions, the theorem of McCarthy and the Busche-Ramanujan identity from the point of view of isobaric representation,   putting these ideas in a different and more transparent light. 

In Section 6,  a structure theorem for the local convolution group of multiplicative arithmetic functions is proved.

In Section 7, we propose a classification system for MF which consists of the categories \textit{degree}, \textit{type}, and \textit{valence}, which enables us to take a refined look at the structure of MF's.  The degree is the degree of the core, the type of an MF has to do with the sizes of its domain and range, and
the valence, with its convolution structure.  We have discussed and extended some results of Laohakosol-Pabhapote, \cite{LP},  Proposition \ref{myproposition8}, Theorem \ref{mytheorem11}  and especially, Theorem \ref{mytheorem13} which gives a candidate for a generalization of the Busche-Ramanujan identity in terms of Schur-hook functions.

In Section 8, we discuss the Kasava Menon Norm for MF in terms of the framework of this paper, showing that it is multiplicative and preserves degree. The proof of multiplicativity has the interesting consequence (Corollary  \ref{mycorollary15}) that it is a special case of a more general fact about convolution products.

In Section 9,  we state the divisible embedding theorem for the WIP-module \cite{MT}, \cite{C} and apply it to the MF's;  that is, we give a formula for producing $q-th$-roots, $q \in\mathbb{Q}$, with respect to the convolution product for each Schur-hook polynomial in the WIP-module.  This construction translates to the convolution group of multiplicative functions,  embedding it in its divisible closure. We also apply some of our results from a previous paper  concerning linear recursions to multiplicative functions, \cite{MW}.
\vspace{0.5cm}

1.  \textbf{RING OF ISOBARIC POLYNOMIALS.}

\vspace{0.50cm}

\begin{definition}\label{mydefinition1}
 An \textbf{isobaric polynomial} is a polynomial of the form $$ \sum_{\alpha \vdash n} C_{\alpha}  t_1^{\alpha_1}... t_k^{\alpha_k},  \quad \alpha = (\alpha_1,...,\alpha_k), \quad \sum j\alpha_j =n, \quad \alpha_j\in \mathbb{N}$$
\end{definition}
\vspace{0.50cm}

The condition $\sum j \alpha = n$  is equivalent to:  $(1^{\alpha_1},...,k^{\alpha_k})$ is a partition of  $n$, whose
 largest part is at most $k$, and written $\mathbf{\alpha} =  (1^{\alpha_1}... k^{\alpha_k})  \vdash  n.$  Thus an isobaric polynomial of \textit{isobaric degree} $n$ is a polynomial whose monomials represent partitions of $n$ with largest part not exceeding $k$. These polynomials form a graded commutative ring with identity under ordinary multiplication and addition of polynomials, graded by isobaric degree. This ring is naturally isomorphic to the ring of symmetric polynomials, where the isomorphism is given by the mapping:

$$t_j \to (-1)^j e_j$$

\noindent with $e_j$ being the j-th elementary symmetric polynomial of isobaric degree  $n$ in the graded ring of symmetric polynomials on the monomial basis \cite{MT}, \cite{Macd} .  This isomorphism associates the Complete Symmetric Polynomials (CSP) in the monomial ring with the Generalized Fibonacci Polynomials (GFP) in the isobaric ring.  It also associates the Power Symmetric Polynomials with the Generalized Lucas Polynomials (GLP) in the isobaric ring \cite{MT, Macd} .  We denote these two sequences of polynomials, respectively,  $$ \{F_{k,n}\} \quad\quad and \quad\quad  \{G_{k,n}\} $$
where $k$ is the number of variables,  and  $n$ is the \textit{isobaric degree}.

They can be explicitly represented as follows:

 $$F_{k,n} = \sum_{\alpha \vdash n} 
  \left( \begin{array} {c}   | \alpha | \\  \alpha_1,...,\alpha_k  \\   \end{array}
\right)  t_1^{\alpha_1}... t_k^{\alpha_k}  . $$
 
 $$G_{k,n} =  \sum_{\alpha \vdash n} 
  \left( \begin{array} {c}  | \alpha | \\  \alpha_1,...,\alpha_k  \\   \end{array}\right) \frac{n} {|\alpha|}  t_1^{\alpha_1}... t_k^{\alpha_k}  . $$
    
 {\large{\Large} For a fixed $k$ both the GFP and the GLP are linearly recursive sequences indexed by   $n$,  and as we shall see below,  the indexing can be extended to the negative integers with  preservation of the linear recursion property \cite{MT, MW}.  Other sequences of isobaric polynomials with the linear recursion property,  consist of elements which are isobaric reflects (i.e., images under the natural isomorphism with the monomial ring) of Schur-hook polynomials (to be explicitly defined below).
It turns out that these sequences can themselves be regarded as  elements of a $Z-$module, called the Module of Weighted Isobaric Polynomials (WIP-module) \cite{MT, MT2}.  It is a remarkable fact that the only isobaric polynomials that can be elements in linearly recursive sequences of weighted isobaric polynomials are those that occur in one of the members (sequences) in the WIP-module (\cite{MT}, Theorem 3.4) . We shall speak of a polynomial as being in the WIP-module if it occurs as a member of a sequence in this module. With this convention we define a WIP-polynomial by the equation

 $$ P_{\omega,k,n} =  \sum_{\alpha \vdash n} 
\left(
\begin{array}{c}
 |\alpha|   \\\alpha_1,...,\alpha_k
\end{array}
\right)
\frac{\sum{\alpha_j\omega_j}}{\sum\alpha_j}  t_1^{\alpha_1}... t_k^{\alpha_k}  . $$
   
 $$ \mathbf{\omega} =  (\omega_1,...,\omega_k)$$
 
 \noindent where $ \mathbf{\omega} =  (\omega_1,...,\omega_k)$ is the $weight$ vector, usually taken to be an integer vector.  Both the GFP's and the GLP's are weighted sequences,  the weighting being given by, respectively, $(1,1,...,1,...)$ for the GFP's,  and $(1,2,...,n,...)$ for the GLP's. The Schur-hook polynomial sequences have weightings of the form $\pm(0,...,1,...,1,...)$.  The GFP, in particular, are Schur-hooks.  The GLP'S
are alternating sums of all of the Schur-hooks of the same isobaric degree. Let  $k$ and $n$ be fixed, then a WIP-sequence is determined by letting $n$ vary in $\mathbb{Z}$.

We need another concept which binds these various sequences together,  namely, that of the \textit{core polynomial}, a concept, though not the name, familiar in the theory of linear recursions.

\begin{definition}\label{mydefinition2}
Given a set of variables $ \mathbf{t} =(t_1,...,t_k), $ the  \noindent \textbf{ \textit{core }} polynomial is the polynomial

\vspace{0.5cm}
$$ (*) \quad\quad[t_1,...,t_k] =  X^k - t_1X^{k-1} - ... - t_k.$$
\end{definition} 

This polynomial is related to the various sequences of isobaric polynomials by the two fundamental theorems of symmetric functions, the first being that the ring of symmetric functions is generated by the elementary symmetric functions,  the second,  that the coefficients of  a monic polynomial are (up to signs)  elementary symmetric functions of the roots.  A rather striking way of immediately achieving this connection is through the $ \mathbf{companion\, matrix}$ (CM).

\vspace{0.5cm}
2. \textbf{THE COMPANION MATRIX AND CORE POLYNOMIALS}

\vspace{0.50cm}
Given the Core Polynomial,  $[t_1,...,t_k] $,  the companion matrix  is the matrix
\vspace{0.50cm}

 \begin{center} $A = \left(\begin{array}{cccc}0 & 1 & ... & 0 \\... & ... &...&1 \\t_k & t_{k-1} & ...&t_1\end{array}\right).$
 \end{center}
\vspace{0.50cm}

A remarkable property of  the Companion Matrix $A$ is that when  $A$ operates on the row-vectors of the matrix,  say on the right, the orbit of $A$ is just the matrix representation of the powers of  $A$. Thus if we list the rows of the orbit of $A$  as an $\infty \times k$-matrix, then every block of $k$-contiguous rows is an $n-th$ power of  $A$.  Moreover, $A$ is invertible iff $t_k \neq 0$,  so that the matrix can be extended northward 
to the negative powers of  $A$, i.e., to $A^{-n}$,whenever $t_k \neq0$, giving a doubly infinite north-south matrix with $k$ columns.  In what follows,  we shall assume that $A$  is non-singular, i.e., $t_k \neq0$. In fact,  it will be convenient to assume that, unless  said otherwise,  the core polynomial is irreducible.  For reasons that will become apparent later, no generality is lost by this assumption.

We can now write the orbit matrix $A^\infty $ described above as follows:

\vspace{0.50cm}

 {$A^\infty = \left(\begin{array}{cccc}...&...&...&... \\(-1^{k-1})S_{(-2,1^{k-1})}&...& -S_{(-2,1)} & S_{(-2)} \\(-1^{k-1})S_{(-1,1^{k-1})} &...&- S_{(-1,1)} & S_{(-1)} \\(-1^{k-1})S_{(0,1^{k-1})}&... & -S_{(0,1)} & S_{(0)} \\(-1^{k-1})S_{(1,1^{[k-1})}&... & -S_{(1,1)} & S_{(1)} \\(-1^{k-1})S_{(2,1^{k-1})} &...& -S_{(2,1)} & S_{(2)} \\(-1^{k-1})S_{(3,1^{k-1})}&... & -S_{(3,1)} & S_{(3)} \\(-1^{k-1})S_{(4,1^{k-1})}&... & -S_{(4,1)} & S_{(4)}\\ ...&...&...\end{array}\right)$.  
 
 \vspace{0.5cm}
 
 This matrix can be regarded as a listing all of the positive and negative powers of the matrix  $A$, \cite{MW},\cite{CL};  as such,  it represents the free abelian group generated by  $A$, an infinite cyclic group.
 
 The entries here are all of the positive and negatively-indexed Schur-hook polynomials,   where the entries can be converted to isobaric polynomials using Jacobi-Trudi matrix representation \cite{Macd, MT}. 
 
 This matrix carries an extraordinary amount of information.  For example, The right hand column is just the sequence $\{F_{k,n}\}$ of GFP's. Any entry gives the Schur-hook polynomial induced by the Young diagram $\pm(n,1^{k-j})$ , the diagram with an arm of length  $n$ and a leg of length $k-j$.  The negatively indexed symbols represent new "Schur-hook" polynomials whose existence is implied by this matrix. Each column is a doubly-infinite k-degree linear recursion.  The sums of the diagonal elements,  that is,  the traces of the elements of the infinite cyclic group generated by $A$,  is just the (k-degree linear recursive) sequence $(\{G_{k,n}\})$, that is, the GLP's.The rows give a representation of the powers of the roots of the core polynomial in terms of a basis consisting of the first $k-1$ powers of a root $\lambda$, namely, $\lambda^0,...,\lambda^{k-1}$.
And the limit of quotients of successive terms of the $F_{k,n}-$sequence can be written in terms of the roots of the core polynomial. (It is well known, e.g., that the quotients of successive terms of the Fibonacci series is either the golden ratio or the negative of the reciprocal of the golden ratio, the other root of the quadratic equation defining the golden ratio.  This is the case where the core polynomial is $[1,1]$ implying that the sequence ${F_n} $ is the Fibonacci sequence.) Similar results hold for the the sequences of Schur-hooks in the other columns of the companion matrix.

One can also ask if the negative Schur-hook polynomials have a description that makes them more transparent,  and also whether there are useful computational tools.  The answer to both questions is yes. 

   \begin{proposition}\label{myproposition1}

Let $k$ be the degree of the core, i.e, we are considering isobaric polynomials on $k$ variables, or equivalently,  a companion matrix with $k$ columns;  and,  let  $s$ be the column coordinate, measuring \textit{from right to left}, with $s=0,1,...,k-1$,  $k-j$, the row coordinate, $k$ fixed,  $j= 0,  1,  2, ...$ .  Then 

$$(-1)^s S_{(-k-j,1^s)} = (-1)^{s+1}\frac{S_{((j+1)^s,j^{k-s-1})}}{t_k^{(j+1)}}$$

\noindent expresses the Schur-hook element of the companion matrix in the $(k-j,s)$ -position in terms of a fraction whose numerator is a positively indexed Schur polynomial, and the denominator, a power of the Schur polynomial $S_{1^k}$. Note that the numerator is a Schur-hook only if $k\leqslant 2$ and $j \leqslant 1$.
  $\square$

\end{proposition} 

Thus we have expressed a negatively indexed element as the quotient of a Schur polynomial by a power of a particularly simple Schur-hook,  the numerator being (except in low degree cases) a non-hook polynomial.   This result will not be used below, so we defer a proof to a later paper; however, using the recursion properties of the sequences,  \textit{or} the isobaric version of the Jacobi-Trudi formula, the reader will have no trouble deriving  the formula.

As for computation,  as we have just pointed out, the columns of this infinite matrix are linearly recursive with recursion coefficients $t_1,...,t_k$; so, in particular, the positive part of a sequence determines the negative part.  And,  each term can be computed individually using the Jacobi-Trudi formula by extending the F-series to negatively indexed terms.   

A natural  question then is how can the negatively indexed Schur-hooks  be interpreted combinatorially. For example, what is the analogue of Young diagram combinatorics for these polynomials,  and what is the conjugate of a negative hook polynomial?  A clue to the answers might lie with Proposition \ref{myproposition1}.

A second observation concerning the negatively indexed elements of the companion matrix is the following:  It follows immediately from the definition of isobaric polynomials that they are polynomials whose monomials are indexed by partitions of an integer $n$ with largest part $k$, where  $n$ is the isobaric degree of the (positively indexed) polynomial.  The rather pleasant fact is that the negatively indexed polynomials in the companion matrix  have monomials indexed by positive and negative partitions of $-n$ with largest part $|k|$.  Multiplying these partitions through by $-1$ gives a version of the \textit{signed partitions} of $n$ as recently discussed by George E. Andrews, \cite{GEW}

We note that as $k$ increases  the identity matrix inside of  $A^\infty$ increases in size, and the non-zero terms of the negative sequences recede. In the limit,  there is no non-trivial negative part to each sequence.

 There are other algebraic structures that can be imposed on the isobaric ring,  \cite{MT}.  For subsequent use, we want to consider one such, a new product on the elements of the WIP-module.  Namely, the convolution product.  

 \begin{proposition}\label{myproposition2}
For a fixed $k$, the\textit{ convolution product} of two elements  $U_n$ and $V_n$ in the WIP-module is the convolution $U_n \ast V_n =  \sum_{j=0}^n U_j V_{n-j}$. 
\end{proposition} 

  It turns out that this gives a graded group structure to the WIP-module.  In particular, $-t_n$  is the convolution inverse of  $F_{k,n}$,  where  $t_j = 0$ when $j > k$  and $t_0 = 1$. The fact that this is the convolution inverse is equivalent to the statement that $F_{k,n} = \sum_{j=1}^k t_jF_{k-j},$ that is, that the GFP-sequence is a $k-th$ order linear recurrence.  

We have regarded the $t_j's$ as indeterminates so far, and the core polynomial as a generic $kth$ degree polynomial.  That  is we have been operating with polynomials, not polynomial functions; but there are many applications in which it is convenient to evaluate these polynomials over a suitable ring. It is in this context that the names Generalized Fibonacci and Generalized Lucas were chosen.  If  $k=2$ and $t_1=1=t_2$,  then the $F-sequence$ is just the Fibonacci sequence,  and the $G-sequence$ is the Lucas Sequence.  For some purposes, we shall want to choose the evaluation ring to be $\mathbb{ Z}$, but other rings will also be useful. 
 \vspace{0.50cm}

We record here the generating functions for elements in the WIP-module. For example,  a generating function for a GFP is given by $$H(y) =  \frac{1}{1-p(y)},  \,\,\; p(y) = t_1y+...+t_k y^k.$$
 where $p(y)$ is the generating function for the convolution inverse of the $F_{k,n}$.  It will be convenient to call the $F_{k,n}$  \textit{positive terms } and the $(t_1,...,t_k), \textit{negative terms}$.
 
For an arbitrary element  $P_{\omega,n}$ in the WIP-module, we have the generating function $\Omega(y) = \sum_{n\geqslant0} P_{\omega,n} y^n$ given in closed form by
  
\begin{proposition}\label{myproposition3}
$$\Omega(y) = 1 +  \frac{\omega_1 t_1 y + \omega_2 t_2 y^2 +...+\omega_kt_ky^k}{1-p(y)}$$

\noindent where $ P_{\omega,0} = 1.$ [MT1]
\end{proposition}
\vspace{0.50cm}

As with the GFP's,  we shall call the sequences $ P_{\omega,n}$ \textit{positive} and their convolution inverses, \textit{negative}.

 In \cite{MW} we regarded the isobaric ring as a ring of functions and studied the sequences in the WIP module with respect to periodicity and periodicity modulo a prime (a linear recursion is periodic if and only if every root of the core polynomial is a root of unity;  and, trivially, every linear recursion is periodic modulo $p$  for every prime p \cite{MW}, theorems 2.1 and 2.2). If $c_p[t_1,...,t_k]$ denotes the period of 
 $F_{k,n}$,  then $c_p[t_1,...,t_k]$ is an invariant of every sequence in the WIP-module, and of the core polynomial. 
(Letting $p=1$ takes care of the case covered by Theorem 2.1). The fact that every sequence in the WIP-module has the same $p-$period has consequences for other structures derived from the same core polynomial. This leads to results concerning the number fields obtained as quotients by an irreducible core polynomial, discussed in \cite{MW}.  Another such application occurs in the ring of Arithmetic Functions,  which we turn to now.
  \vspace{0.5cm}

 3. \textbf{THE RING OF ARITHMETIC FUNCTIONS}
\vspace{0.50cm} 

While the isobaric ring is a not so classical version of the well known ring of symmetric functions,  the elements in the ring of arithmetic functions have long been objects of study, but not usually from a structural point of view (but see \cite{CE},\cite{R},\cite{R2}\cite{HS}, and recently, \cite{LP}, \cite{LPW}, \cite{LPW2},\cite{PH},\cite{PH2},\cite{PH3},\cite{PH4},\cite{PH5}).  It is possible that the relation between the two structures is implicitly well understood,  but it is rather surprising that the relationship, to our knowledge, has not been made explicit in the literature.  But see \cite{TM} , \cite{TM2} and \cite{MT}. The connection is that the GFP's with the convolution product  is locally isomorphic to the group of  multiplicative functions  under the convolution product,  and,  by consequence, every sequence in the WIP-module yields a group of multiplicative functions that can be associated locally with a given multiplicative function.  

So now we shall review the facts about arithmetic functions that we need in order to show this connection. 
We recall that  arithmetic functions $(\mathcal{A})$ are functions from $\mathbb{N}$ to $\mathbb{C}$,  and form a ring under the usual sum and product rule for functions;   but it is also usual to add convolution as an additional operation:  $(\alpha \ast \beta) (n) = \sum_{d|n} \alpha (d) \beta(n/d)$.  Call this structure in which the convolution product substitutes for the standard product of functions, $\mathcal{A}^*$. Then  $\mathcal{A}^*$ is also a ring,  a commutative ring with identity, and, in fact, a unique factorization domain \cite{CE}.  An arithmetic function $\alpha$ is invertible with respect to convolution iff  $\alpha(1) \neq 0$.  $\alpha$ is a\textbf{ multiplicative function} (MF) iff $\alpha(mn) = \alpha(m) \ast \alpha(n)$ whenever $(m,n) = 1$.  Thus the MF's are just those AF's that are uniquely determined at powers of primes.  It follows from all of this,  that if an MF is invertible, its value at $1$ is $1$.  The zero-function is an arithmetic function.  We shall exclude it from the set  MF, and note that, with this exclusion,  MF belongs to the group of units of   $\mathcal{A}^*$ \cite{PJM}.

The multiplicative functions which are units in   $\mathcal{A}^*$ are well known to constitute a group.  Call it  $\mathcal{M}$. An MF $\alpha$ is $\mathbf{completely \;multiplicative}$ if  $\alpha(mn) = \alpha(m) \alpha(n)$ for all $m,n \in\mathbb{N}$ . The subgroup of the group of units in  $\mathcal{A}^*$ generated by the 
 completely multiplicative functions, call it  $\mathcal{L}$ , is known to be a(n) (uncountably generated) free abelian group \cite{TM}.  $\mathcal{L}$ is clearly a subgroup of  $\mathcal{M}$.  If  $\alpha \in \mathcal{L}$, then inverses are determined at each prime,  that is, $\textit{locally}$, either  by a monic polynomial of degree $k$, where $k$ depends on the prime $p$, or by a polynomial of infinite degree, that is, by a power series. Call either of these the \textit{local core} of the MF, denoted by $ \textsl{C}_p(X) $.  The restriction of the invertible functions in $\mathcal{M}$ to their values on a particular prime $p$ form a subgroup of $\mathcal{M}$, which we shall denote $\mathcal{M}_p$.  In what follows we shall often drop the subscript $p$ with the understanding that our discussion is local. We do this in the following definitions. A multiplicative function is called \textit{positive}  if it is the convolution product of CM functions, \textit{negative} if it is the product of the inverses of CM  functions, and it is called \textit{mixed} if it is the convolution product of at least one non-identity CM function and one negative of a non-identity CM function.  We can also classify functions in MF as being of one of the four following types depending on the sizes of their ranges and domains: 
\newpage \begin{enumerate}
\item  type (1)   $(fin,fin)$;
\item  type (2)   $(\infty, fin)$ ;
 \item  type (3)   $(fin, \infty )$;
 \item  type (4)   $(\infty,\infty).$
\end{enumerate}

The notation here means the pair $(range, domain)$ = $(F\textbf{(t)},\textbf{t})$.
The first type are those MF's which have both range and domain finite; the second type, those that have an infinite range  and a finite domain; the third type, those with a finite range and an infinite domain; and the fourth type has both range and domain infinite--- \textit{ finite range} and \textit{finite domain} mean that eventually all values are zero.  In Theorem \ref{mytheorem10} and Corollary \ref{mycorollary12} we show that the set of type (1) functions contains only the identity, that type (2) functions are just the \textit{positive} MF's, that type (3), the negative MF's; and that type(4) are \textit{mixed}. We shall discuss this point further in section 7.
  \vspace{0.5cm}
 
 4.\textbf{ RELATION BETWEEN MULTIPLICATIVE FUNCTIONS AND THE WIP-MODULE}
   \vspace{0.5cm}
 
 Since a given core polynomial determines both  the WIP-Module (in particular,  the infinite companion matrix),  and determines a particular arithmetic function locally in  $\mathcal{M}$, it is clear that there is a strong connection between the MF's and the WIP sequences.  In fact, the GFP-sequence evaluated at the vector $\textbf{t}$ is a (non-trivial) MF. So is the positively indexed part of every column in the matrix $A^\infty$.  Thus every sequence in the WIP-module is by consequence also in MF. There are instances of MF's for which the core polynomial is constant for all choices of prime (the MF  $\tau$, of degree 2 which counts the number of divisors of $n$ has the core polynomial $X^2-2X+1$). There are also instances where the core polynomial has the same degree over all primes and the coefficients are given by the same functions of $p$ (the MF $\sigma$ of degree 2 which records the sum of the divisors of $n$ where the core polynomial is given by $X^2-(p+1)X+p$.  The first case can be regarded as a special case of the second. It prompts us to consider the  \textbf{generic }core polynomial $X^k+\sum_{j=1}^k - t_jX^{k-j}$.  
 First we note that  convolution preserves isobaric degree in the WIP-module, and core degree in both the isobaric case and the $\mathcal{L}$ case. The analogue of \textit{ isobaric degree} for the multiplicative functions is just \textit{power of the prime} $p$.  The analogue of a function in $\mathcal{L}$ requiring  infinitely many powers of the prime in its definition is that $k$ is unbounded in the WIP-module; that is that  ($t_j$) is different from $0$ for infinitely many  $j's$, or, equivalently,  that the companion matrix is unbounded on the west.  If we  call the MF's which are locally of degree $k$ for all primes $p$, $k$-\textit{uniform}, and the set of all $k$-uniform for all $k$, \textit{uniform}, then it is easily seen that the uniform MF's form a  graded group under convolution.  At each level of the grading and for each prime $p$ the core polynomial induces a (cyclic ) direct summand (the values of the MF at the powers of the prime $p$)  and, at the same time, on the induced GFP.  That is, for each $k$ the subgroup at that level is a direct sum of cyclic subgroups, one such subgroup for each  $p$.  We refer to these subgroups as the local subgroups of degree $k$. These subgroups all have the same generic core, i.e., the cores of the elements in the subgroup are evaluations of the same generic polynomial. So far  we have glossed over the point that the output of the functions $F_{k,n}$ depends upon the choice of domain ring.  It is clear that the output of such evaluations will always be elements in the same ring as that of the input; for examples an input of integers yields an output of integers. Integer inputs will often be used in the examples simply because so many classical MF's are of that sort. We shall stretch the term \textit{counting function} to include any multiplicative function with an integer output.  In general, if the evaluation ring is $\mathcal{R}$  (a  subring of the complex numbers), we denote the subgroup in $\mathcal{M}$ that they generate,  $\mathcal{M_R}$.  
 {0.5cm}e{0.5cm}
 
 With these remarks, we state the main theorem.

\begin{theorem}\label{mytheorem1}
Given a prime $p$, let $ \chi \in \mathcal{M}$ with local-core $\textsl{C}(X), \, k $ finite or infinite, and let $F_{k}(\mathbf{t})$ be the sequence of GFP's induced by this core,  then  $F_{k,n}(\mathbf{t}) = \chi(p^n) $.
\end{theorem}

\textbf{Proof}   As pointed out above, for each integer $k$, every polynomial in the WIP-module determined by $k$  (i.e., whose companion matrix is $\infty \times k$) is a member of a k-linear recursive sequence;  in particular,  the GFP-sequence is one such sequence. Thus, $F_{k,n} = \sum_{j=0}^k t_jF_{k,k-j}, F_{k,0} = 1$, where $t_j, j=1,...,k$ is the set of parameters that determine the recursion. Thus , every linear recursion of degree $k$ is determined by choosing a set of values for $t_j$. \cite{MW}.  ($k$ can be finite or infinite).  It is clear such a choice of parameters determines a multiplicative arithmetic function locally---the recursive relation determines the convolution product. 

 The converse is also true; that is, each multiplicative function  $\chi$ has a locally faithful representation as an evaluation of $F(t_1,...,t_j,...)$ in the GFP-sequence.  For, given a prime $p$ and the set of values  $\chi(p^n) =a_n$,  we can determined the $t_j$ and $F_{k,j}(t_1,...,t_j,...)$ inductively  in such a way that  $F_{k.j}(t_1,...,t_j)=  \chi(p^j), j\leqslant k$. 
 Let $ \chi(p^j) =a_j, j = 1,2,...$, and let $F_{k,0}=a_0 =1 $.  Let $t_1 = a_1 = F_{1,1}$. Suppose that $t_j$ for $j<n+1$ has been defined and that $a_j = F_{k,j} j<n+1$.  We define $t_{n+1}$ by  
$$t_{n+1} = a_n -\sum_{j=1}^n t_j a_{j-1}.$$  

That is, $$t_{n+1} = \chi(p^n) -\sum_{j=1}^{n-1}t_j \chi(p^{j-1})$$ (cf. Proposition \ref{myproposition5})
$$ = F_{n+1,n+1} -\sum_{j=1}^{n}t_j F_{n,n-j+1}.$$

The theorem now follows by the recursive property of the GFP sequence and induction.  $\square$

The last equation in the proof reflects the fact that $F_{k+1,k+1} - F_{k,k+1} = t_{k+1}$. We refer to this property of the GFP sequence as a \textit{conservation} principle.  The following list containing the first few terms of the GFP is useful for understanding these remarks:\begin{itemize}
  \item $F_{k,1} = t_1$
  \item $F_{k,2} = t_1^2 + t_2$
  \item $F_{k.3} = t_1^3 + 2t_1t_2+t_3$
  \item $F_{k,4} = t_1^4 + 3t_1^2t_2 + t_2^2+2t_1t_3+t_4$
  \item $F_{k,5} = t_1^5 + 4t_1^3t_2 +  3t_1t_2^2 +3t_1^2t_3+2 t_2t_3+2t_1t_4+t_5$
\end{itemize}

\noindent if $k\leqslant 4$.  One can compare these GFP polynomials with the values determined by the closed formula in Section 1.

The very fact that the construction in the proof of Theorem \ref{mytheorem1} is possible guarantees, by the way, that every multiplicative function is recursive.

\begin{definition}\label{mydefinition3}
Let $ \mathcal{M}_p$ be the subgroup of  (the graded group) $\mathcal{M}$ of multiplicative functions restricted to the prime $p$.
\end{definition}
Consider the set of all arithmetic functions determined locally by evaluating the set of isobaric polynomials $\{F_k(\textbf{t})\}$ over the complex numbers. This is a graded group which we denote by  
 $\widetilde{\mathcal{M}}_p$.
 
 \begin{corollary}\label{mycorollary2}
 $\widetilde{\mathcal{M}}_p$  is isomorphic to $\mathcal{M}_p$.       $ \square$
\end{corollary}

\vspace{0.5cm}
 It is useful to record generating functions for the $F-$sequences, and hence for the associated MF's, as well as for arbitrary polynomials in WIP. 
  \begin{theorem}\label{mytheorem3}
  
Let $X^k-t_1X^{k-1}-...-t_k = \textsl{C(X)}$ be a generic core,  and let $P(y) = t_1y +...+t_ky^k$,  then 
$$H(y) = \frac{1}{1-P(y)}$$ is a generating function for the sequence $\{F_{(k,n)}(\textbf{t})\}$,  And, in general, for any element of WIP $$W(y) = 1+ \frac{\omega_1t_1y + ... + \omega_kt_ky^k}{1-P(y)}, $$where \textbf{$\omega$} is the weight vector for the element of  WIP.  (\cite{MT}, Theorem 4.2, Lemma 4.3)
\end{theorem}

\begin{remark}\label{myremark1a}

If we denote the group of units of $\mathcal{A}^\ast$ by $\mathcal{U}$,  then clearly  $\mathcal{L} \leqslant \mathcal{M} \leqslant \mathcal{U}$. The question of when the relation is that of equality arises.
In the case of $\mathcal{L}$ and $\mathcal{M}$, it is clear that  (non-zero) multiplicative functions
of types 1,2, and 3 are in  $\mathcal{L}$. It will be clear in what follows that  $\mathcal{L} \cap type \,4 \neq \emptyset$ (Remark \ref{myremark2}, Proposition \ref{myproposition7}). We conjecture that if $\alpha \in \mathcal{M}$ and is of type 4, that $\alpha\in \mathcal{L}$. The truth of this conjecture would imply equality of  $\mathcal{L}$ and  $\mathcal{M}$.  We would also conjecture the equality of  $\mathcal{M}$ and  $\mathcal{U}$, but this is probably a more difficult problem.
\end{remark}
 
5. \textbf{EXAMPLES}
 \vspace{0.5cm}

Consider the multiplicative functions $\tau$ and $\sigma$, where $\tau$ is the function on $\mathbb{N}$
which counts its distinct divisors,  while $\sigma$ is the divisor sum function.  Both of these 
functions are multiplicative of degree  2. $\tau(p^n) = n+1$ and   $\sigma(p^n) = 1+p+...+p^n$.  We can find the local core polynomial for $\tau$ by noting that $\tau(p) = 2$, $\tau(p^2)= 3$ and $\tau(p^3) = 4$.  Using the induced GFP, $F_1 = t_1$,  we have that $t_1 = \tau(p^1) = 2$. Then we note that $F_2 = t_1^2+t_2 = \tau(p^2)= 3$,  we deduce that  $t_2 = -1$. A similar computation for $t_3$ yields $t_3=0$. An induction using the recursive properties of the GFP sequence, shows that $t_j=0$ for all $j>2$.  Thus the local core is the quadratic polynomial  $X^2-2X+1$,  which, incidentally shows the well known result that  $\tau$ is the convolution product of two copies of $\zeta$,  all of whose values are $1$.

If we carry out the same procedure for  $\sigma$, we find that $t_1 = 1+p$, that $t_2 = -p$ and that the degree of $\sigma$ is $2$, that is, that $t_j=0$ for $j>2$, hence, the local core is 
$X^2-(p+1)X+p$. Again, since the local core has linear factors $X-p, X-1$, $\sigma$ is the convolution product  $\zeta_1 \ast \zeta$ of two
CM  arithmetic functions,  i.e., two degree $1$  functions,, where    $\zeta_k(p^n)  = p^{nk}$.
Degree 2 uniform MF's are also called \textit{specially multiplicative}.These are both examples of uniform MF's, i.e., they are elements of $\mathcal{L}$.  This same procedure applied to the Euler totient function $\phi$ shows that $t_j\neq 0$ for all $j>0$.  Thus $\phi$ is  uniform  and is an example of a MF whose core is a power series.  Its values are given by $F_{k,n}(t_1,...,t_k,...), t_j=p-1$ for all $j>0$ and all $k > 0$ and all primes $p$.  It is well known that $\phi = \zeta_1 \ast \mu$, where $\mu$ is the convolution inverse of $\zeta$, i.e.  $\phi = \zeta_1 \ast \zeta^{-1}$, which is called in \cite{LP} a function of valence $<1,1>$ (see Definition \ref{mydefinition4} below), part of a general notation $<r,s>)$ for multiplicative functions which are convolution products of $r$ completely multiplicative functions and $s$ inverses of completely multiplicative functions. We shall have more to say about such functions later, but prefer to change the notation so as not to conflict with notation already established here. 

\begin{definition}\label{mydefinition4}
A multiplicative function has \textit{valence} $<r,s>$ if it is a convolution product of $r$ completely multiplicative functions and $s$ inverses of completely multiplicative functions.
\end{definition}
So in this notation $\phi$ is of valence $<1,1>.$ This leads to an interesting theorem.

\begin{theorem}\label{mytheorem4}
Let $\alpha = \beta \ast \gamma^{-1}$ where $ \beta$ and $\gamma$ are CM,  that is,  $\alpha$ is of valence $<1,1>$,  then $\alpha$ is of type $(inf , inf)$.  That is $\alpha$ has an infinite range and infinite domain;  thus $\alpha$ has an infinite core.
\end{theorem}

\textbf{Proof}  We combine the techniques of calculation used before with the generating functions from Theorem \ref{mytheorem3}. Represent $\alpha, \beta, \gamma$ by $F, F', F''$, respectively.  Letting the parameters for the two CM functions $\beta, \gamma$ be given as $t'_1 \neq 0, t_j = 0$ otherwise, and  $t''_1\neq 0, t_j = 0$ otherwise;  then the local core polynomials are $X-t'_1, X-t''_1$.  From  Theorem \ref{mytheorem3} we have that the generating function for the convolution product in terms of the parameters of the factors is 
$$\frac{1-t''_1y}{1-t'_1y} $$ 
$$= 1+ \sum_{n=1}^{\infty} ( t'^n_1 - t'^{n-1}_1t''_1)y^n$$
$$=\sum_{n=0}^{\infty} F_ny^n.$$
Thus $$F_n = ( t'^n_1 - t'^{n-1}_1t''_1).$$

Using these values for $F_n$, the calculating methods employed above, and induction, we can deduce that $$t_n = (- t''^n_1 + t'_1t''^{n-1}_1).$$

But, since $t'_1 \neq 0$  and  $t''_1\neq 0$,  this shows that both the range and domain of $F_n$ are infinite, that is,  $\alpha$ is of type $(inf , inf)$. $\square$

Thus the case for $\phi$ generalizes.  It seems reasonable to conjecture that the correct generalization is to $<r,r>)$-functions.

  It is also instructive to look at the convolution product  $\tau\ast \sigma = \alpha $, a positive function.  A calculation of the nature of those above shows that the parameters for  $\tau \ast \sigma = \alpha$ are :  $u_1 = p+3, u_2 = -3(p+1), u_3 = 4p+1, u_4 = p$, all non-zero, and $u_j= 0, j>4$.  Hence,  the product is of degree 4,  and the local core is $X^4-(p+3)X^3-3(p+1)X^2-(4p+1)X-p$, which is just the product of the local cores of the two factors.  We have then that 

\begin{definition}\label{mydefinition5}
For any MF,  we define its \textit{degree} to be $k$ if the parameters $t_j$ of the function have the property that $t_k\neq 0$ and $t_j = 0, j>k$, where k is either finite or infinite.
\end{definition}

The previous example prompts the following two theorems
 
 \begin{theorem}\label{mytheorem5}
 $$   core(\alpha_1 \ast \alpha_2)=core(\alpha_1) core(\alpha_2) \quad\quad \square $$
\end{theorem}
\begin{theorem}\label{mytheorem6}
 $$deg(\alpha_1 \ast \alpha_2)= deg(\alpha_1) + deg(\alpha_2)  \quad\quad \square$$
\end{theorem}

 The proofs of these two theorems follow easily from generating function considerations. 

 \begin{corollary}\label{mycorollary7}
The positive part of each column in the infinite companion matrix is also an MF,  hence in $\mathcal{M}$. $\square$
\end{corollary} 

 The negative part of each column in the infinite Companion Matrix  is also in MF.  So we have the following situation.  The core polynomial determines, what might be called, a principal MF,  the one determined in Theorem \ref{mytheorem1} by the GFP,  and at the same time what might be called (a module of) satellite MF's, those determined by the other sequences in the WIP-module.  In addition to this,  this connection also provides each such positively indexed sequence of MF's with a negatively indexed sequence of arithmetic functions.  How is one to interpret these negatively indexed arithmetic functions, the analogous question to the one asked above,  how is one to interpret the combinatorics of the negatively indexed Schur-hook polynomials ? 

  The fact that each column in the infinite companion matrix is a k-linear recursion is reflected in the fact that the associate MF's are determined locally by the first $k$ powers of the prime, while the rest of the sequence is determined by linear recursion, the recursion constants being given by the vector $\mathbf{t}$.
    In \cite{MW} it was shown that, in the case that the local core is irreducible,  there is a strong relation between the WIP-Module and the number fields associated with the field extension determined by the local core.  This fact gives a three-way relation among the three structures:  WIP-Module, multiplicative arithmetic functions and number fields. In particular,  it associates with every such number field, a special set of multiplicative functions.

Another question arises from the fact that the UFD $\mathcal{A}^*$ has a rich ideal structure \cite{C,HS}.  Is there a representation of these ideals in terms of symmetric polynomials?  
 \vspace{0.5cm}
 
 6. \textbf{SPECIALLY MULTIPLICATIVE ARITHMETIC FUNCTIONS}
\vspace{0.5cm}

A theorem due to P.J. McCarthy \cite{PJM} states that  a multiplicative function $\chi$ is specially multiplicative, that is, is of degree 2, if and only if for each prime $p$, $$\chi(p^{n+1}) =  \chi(p)  \chi(p^n) - \chi(p^{n-1}) B(p),$$  where $B (p)= \chi(p)^2-\chi(p^2)$, and $B(p) \in CM$. Furthermore,  degree 2 multiplicative arithmetic functions are characterized by the property that they admit a Busche-Ramanujan identity \cite{KGR}. Using Theorem 1 to translate these results into isobaric form,  we get as a characterization of degree 2 MF's the following:   
$$  F_{2,n+1}(t_1,t_2) = t_1 F_{2,n}(t_1,t_2)+t_2 F_{2,n-1}(t_1,t_2),$$

 \noindent or more succinctly,
 
 $$  F_{2,n+1} = t_1 F_{2,n}+t_2 F_{2,n-1} $$

 \noindent or
  
  $$ F_{n+1} = t_1 F_n+t_2 F_{n-1}.$$

\vspace{0.25cm}
In particular,   $B_{\chi} = B (p)=\chi(p)^2-\chi(p^2)$ translates into $-t_2$. But (from the point of view of this paper)  this is just the redundant statement  that degree 2 cores induce linear recursions of degree 2. 

It is also asserted in the McCarthy theorem that $B(p)$ is a completely multiplicative arithmetic function, that is,  that it has degree 1.  Suppose that we represent the degree 2 MF $\chi$ by the GFP-sequence $F(t_1,t_2)$ and $B$ by $\tilde{F}(u_1,...,u_j)$.

The claim concerning $B$ is equivalent to the claim that the parameters $\{u_j\}$ are $u_1 \neq 0$ and $\;u_j = 0, j>1$.  Noting that $B(p) = \chi^2(p)-\chi(p^2) = -t_2 = F_1^2-F_2$ and hence that $B(p^n)  =    \chi(p^n)^2-\chi(p^{2n}) = F_n^2 + F_{2n}$.  Since $\chi$ has degree 2, we can represent it by $\chi = \beta\ast \gamma$, representing $\beta \; and \; \gamma$ by $F'(t'_1) \; and \; F''(t''_1)$.  It then follows that $t_1 = t'_1+t''_1$ and $t_2 = t'_1t''_1$.  Thus $u_1 = -t_2 =t'_1t''_1$. $\tilde{F}_1 = u_1, \tilde{F}_2 = F_2^2-F_4 = (t'^2_1 +  2t'_1t''_1 + t''^2_1)(t'_1t''_1) = u_1^2+u_2 $. From this we can calculate $u_2$. Proposition  $u_2 = (t'_1t''_1)(t'^2_1 + t'_1t''_1+ t'^2_1) = -t_2F_2$. 

The function $\tau$ that counts the number of divisors of $n$ is well known to be a degree 2 function, i.e. specially multiplicative.  If $\chi = \tau$, then $t_1 = 2$ , $t_2 =-1$ and $F_2 = 3$, thus $u_2 = 3$, which gives a counter-example. In fact, we can see that B is deg 1 only if $F_2 =0$.

The Busche-Ramanujan identities for the specially multiplicative functions $\sigma_k$ are the two statements:

$$(1) \sigma_k (mn) = \sum_{d| (m.n)} \sigma_k (\frac{m}{d})\sigma_k \frac{n}{d} \mu(d)d^k$$ 

$$(2)\sigma_k(m) \sigma_k(n) =  \sum_{d| (m.n)} d^k\sigma_k (\frac{mn}{d^2}) $$

where $\sigma_k (n) = \sum_{d|n}d^k$,   $\sigma_1 = \sigma$,  $ \sigma_k = \zeta_k \ast \zeta$.

\vspace{0.5cm}

We translate this into isobaric notation. For simplicity,  we take the case where $k=1$, that is,  $\sigma_1 =  \sigma$.  So we are interested in the identities:

$$(1') \sigma (mn) = \sum_{d| (m.n)} \sigma (\frac{m}{d})\sigma (\frac{n}{d}) \mu(d)$$ 

$$(2')\sigma(m) \sigma(n) =  \sum_{d| (m.n)} d\sigma (\frac{mn}{d^2}) $$
Letting $t_1 = 1+p,  t_2 = -p$, the core coefficients for  $\sigma$, and using Theorem \ref{mytheorem1}, that is
that $F_n(t_1,t_2) = \sigma(p^n)$, since $\sigma$ is a MF of degree 2, and letting $m = p^r, n = p^s, r\leqslant s$,  the two identities become 

\begin{proposition}\label{myproposition4}(Busche-Ramanujan) \cite{PJM}, \cite{RS}                  

$$(4.1) \;F_{2,r+s} =  F_{2,r}F_{2,s} + t_2F_{2,r-1}F_{2,s-1} $$

$$(4.2)\; F_{2,r}F_{2,s}  = F_{2,r+s}  -t_2F_{2,r+s-2} + ...+(-t_2)^j F_{2,r+s-2j}+...+(-t_2)^r F_{2,s-r}$$

where $k$ is the degree of the core, in this case $2$, and $j = 1,...,r$.

\end{proposition}

It is perhaps instructive to give a proof of these well known relations in terms of the isobaric representation of MF's being discussed in this paper.

\vspace{0.5cm}
\textbf{Proof}  We consider the first of these identities.  Omitting the degree index  $k = 2$ on the GFP-symbols, and noting that (4.1) is true when $r=1,$ we have the basis for an induction. But 
\begin{center}
$$F_{r+s}= t_1F_{r+s-1} + t_2    F_{r+s-2}$$ $$=t_1F_{(r-1)+s} + t_2    F_{(r-1)+(s-1)}$$
$$= t_1(F_{r-1}F_s + t_2F_{r-2} F_{s-1}) +t_2F_{r-1}F_{s-1} + t_2^2 F_{r-2} F_{s-2} $$
$$=t_2 (t_1 F_{r-2}F_{s-1} + t_2 F_{r-2}F_{s-2}) + t_1F_{r-1} F_s + t_2 F_{r-1} F_{s-1}$$
$$=t_2(F_{r-2}F_s) +t_1 F_{r-1}F_s +t_2 F_{r-1}F_{s-1}$$
$$= F_s (t_1F_{r-1} + t_2F_{r-2} +t_2F_{r-1}F_{s-1})$$
$$=F_sF_r + t_2F_{r-1}F_{s-1},$$
\end{center}

\vspace{0.25cm}
\noindent using only linear recursion and the induction hypothesis. This proves the first of the two identities. $\square$

\begin{lemma}\label{mylemma1}
$$t_2F_{2,r-1}F_{2,s-1} + \sum_{j=1}^r (-t_2)^j F_{2,r+s-2j} = 0, \;\;  r\leqslant s.$$
\end{lemma}
\vspace{0.25cm}

\textbf{Proof}  We observe that the identity of the lemma holds when  $<r,s> = <1,1>$. Suppose that it also 
holds for $ 2<r+s < n$.  Then, using the linear recursion property of the $F-$sequence, we have, for $k=2$

$$t_2(t_1F_{2,r-1}F_{2,s-2} +t_2F_{2,r-1}F_{2,s-3} )+ \sum_{j=1}^r((-t_2)^j(t_1 F_{2,r+s-2j-1}+(-t_2)^j F_{2,r+s-2j-2} )$$

$$ =(t_2(t_1F_{2,r-1}F_{2,s-2}) +\sum_{j=2}^r(-t_2)^jt_1 F_{2,r+s-2j-1})+(t_2(t_2F_{2,r-1}F_{2,s-3} )+\sum_{j=1}^r (-t_2)^j F_{2,r+s-2j-2} )= 0$$  $\square$

Combining  (4.1) with the Lemma gives   

\vspace{0.25cm}
\begin{corollary}\label{mycorollary8}
(4.1) is equivalent to(4.2). $\square$
\end{corollary}

\vspace{0.05cm} 

It is clear from the proofs of the proposition and  Lemma \ref{mylemma1} that the only assumption made was that we are dealing with degree 2 multiplicative functions, so the results do indeed hold for all degree 2 MF's;  moreover, it is hardly surprising that they do not hold for higher degree MF's,  since our very assumption is that our linear recursions and hence our core is of degree 2, that $t_j = 0, j>2.$ (also see \cite{RS2}, p.282). 
 
 Can the McCarthy characterization of specially multiplicative functions,  that is,  degree two functions, be generalized to finite higher degree functions?  The B-function in McCarthy's theorem is just $-t_2$ at $p$, and the relation itself is just a statement of the linear recursive property of the $F$-functions that represent the multiplicative functions (see (4.5) and ff.).  If we think of $t_2$ as the isobaric degree 2 term in the specially multiplicative case,  then it is natural to think of $t_k$ as the isobaric degree $k$ term in the general case where the core polynomial has degree $k$.  The analogue to the fact that $B(p) = F_1^2 - F_2 = -t_2$ is the following proposition which expresses the indeterminates, $t_j$ in terms of the GFP polynomials, $F_{k,n}$.  
\begin{proposition} \label{myproposition5}
$t_n = F_n(F_1,..., (-1)^{j+1}F_j,..., (-1)^{n+1}F_n),\;\; j=1,2,...,n.$
\end{proposition}

  Thus, \begin{enumerate}
  \item  $t_1 = F_1,$
  \item  $-t_2 = F_1^2-F_2,$
  \item  $t_3 = F_1^3-2F_1F_2+F_3,$
  \item  $...$
\end{enumerate} 
where it is understood that the absence of the subscript  $k$ means that $k$ is allowed to go to $\infty$.
In words,  the  proposition says, write $F_n$ as a function of the $t_j$ then make the substitution  $t_j =  (-1)^{j+1}F_j$.  The construction used here is just the one used in the proof of Theorem \ref{mytheorem1} $\square$

\begin{remark}\label{myremark1}
The main point here is that there is nothing special about specially multiplicative functions.
Expressions of the sort  contained in the McCarthy theorem are true for every multiplicative function and merely represent the fact that these functions are degree k recursive.  The $B(p)-$term in the degree 2 case (that is, $-t_2$) can be replaced by  $-t_2,...,-t_k$ in the degree $k-$case. The representation of $B$ in terms of the original function generalizes to the representation $t_2,...,t_k$ as a function of its associated F-polynomials as indicated in Proposition 5 above. 
\end{remark}  
Proposition \ref{myproposition4} expresses the Busche-Ramanujan identities in terms of the GFP-representation, which in turn suggests a way of generalizing such identities to MF's of higher degree.  Thus,
one way to think of the Busche-Ramanujan identities is as an expression of  $F_{r+s}$ in terms of $F_r$ and $F_s$ together with a \textit{remainder} term. In Theorem \ref{mytheorem13} we have just such a generalization, which has the pleasant property of involving Schur-hook functions as coefficients.

  \vspace{0.5cm}

7. \textbf{STRUCTURE OF THE  CONVOLUTION GROUP  $\mathcal{M}$  OF MULTIPLICATIVE FUNCTIONS REVISITED} 
  \vspace{0.5cm}
  
The group  $\mathcal{L}$ generated by the completely multiplicative functions, sometimes called the group of rational functions, \cite{JR}, \cite{LP}, contains elements of four kinds: the identity, positive elements (the semi-group generated by CM functions), negative functions (the inverses of the postive functions), and mixed elements (those which are convolution products of both positive and negative elements),  as discussed in Section 2.  Each element has a degree which is either infinite or a non-negative integer.  The identity has degree $0$,  a positive element  has positive degree, the degree of its core polynomial, or, equivalently,  the number of CM functions of which it is a product (in \cite{TM} it is shown that the CM functions freely generate a free abelian group).  Both  negative functions and  mixed functions have infinite degrees.  Negative functions have power series cores, and a mixed function has a rational function for a core whose numerator is the core of the positive part and whose denominator is the core of the negative part.  

In Section 2,  the classification of elements of types 1 $(fin,fin)$; 2, $(\infty,fin)$;  3, $(fin,\infty);  $4,$(\infty,\infty)$ in the group $\mathcal{M}$ was introduced,  where  $finite$ range or domain means that the sequence is eventually constantly zero. Infinite means that the sequence has infinitely many non-zero values. 
Clearly, the types are  mutually exclusive; moreover, each type is non-empty. For example, the identity function is type 1. Type 2 consists of the positive functions, and type 3,  the negative functions. Completely multiplicative functions, e.g., $\zeta$, or any specially multiplicative function is of type 2, e.g., $\sigma$, while $\mu$, a negative function, is of type 3,   And, as we shall see, the Euler totient function, $\phi$, is of type 4 .  

  We shall say that the set  of indeterminates  $\{t_j\}$ is the set of parameters of $F$, having in mind that when we assign values to these indeterminates, we determine a particular numerical sequence $F(\textbf{t})$, and, hence, a particular MF.
 
\begin{proposition}\label{myproposition6}
 Let $\alpha \in MF $ and let $\alpha \leftrightarrow F $ and let $\{t_j\}$  be the set of parameters for $F$.  Let $\{s_j\}$ be the set of parameters for $F^{-1}$ which represents $\alpha^{-1}$, then $$\forall j,  \; F_j = -s_j, $$ and $$\forall j,  \; F_j^{-1} = -t_j.$$
\end{proposition}

 \textbf{Proof}  $0 = F_1 \ast  F_1^{-1} = F_1 + F_1^{-1}  = t_1 + s_1$, so $ t_1 = -s_1$;  therefore,
 $$ F_1 = -s_1, \; and \; F_1^{-1} = -t_1.$$
 So assume inductively that $F_j = -s_j, j = 2,...,n-1$ and that $F_j^{-1} = -t_j, j = 2,...,n-1$, then we have that  $$F_n = \sum_{j=1}^{n-1} t_jF_{n-j} +t_n $$  $$ = -\sum_{j=1}^{n-1} t_js_{n-j} +t_n .$$
 
 Similarly,
 $$F_n^{-1} = -\sum_{j=1}^{n-1} t_js_{n-j} +s_n, $$
 therefore, $$F_n - F_n^{-1} = t_n-s_n .$$
 
 Also, 

 $$F_n \ast F_n^{-1} =  F_n + \sum_{j=1}^{n-1} F_{n-j}F_j^{-1}  + F_n^{-1}$$
 
 $$=  F_n - \sum_{j=1}^{n-1} s_{n-j}t_j^{-1}  + F_n^{-1}  = 0.$$
Thus $$ F_n + F_n^{-1}  = - \sum_{j=1}^{n-1} t_js_{n-j} = F_n - t_n =  F_n^{-1} - s_n.$$
 
\noindent So $$2F_n = F_n - t_n $$
 implies  $$F_n = -s_n ,$$
 
 \noindent and similarly, $$F_n^{-1} = - t_n. \quad\quad \square$$  

\begin{corollary}\label{mycorollary9}

 Let $F = F' \ast F''$ with $t'_j, t''_j$ being the parameters of $F', F''$ and $s'_j,s''_j$ the values of $F', F''$, then
$$-s_n = -s'_n +\sum_{j=1}^{n-1}s'_{n-j}s''_j-s''_n$$
$$-t_n = -t'_n +\sum_{j=1}^{n-1}t'_{n-j}t''_j-t''_n$$
equivalently
$$F_n = -s'_n + \sum_{j=1}^{n-1} s'_{n-j}s''_j -s''_n;$$ and  
$$t_n =  t'_n - \sum_{j=1}^{n-1} t'_{n-j}t''_j +t''_n.$$

\end{corollary}

\noindent\textbf{Proof}

Expand the convolution product  $F_n =  \sum_{j=0}^n F'_{n-j}F''_j$ and apply the theorem to the factors. $\square$
 
\begin{definition}\label{mydefinition6}
A convolution product of $r+s$ factors is said to be in \textit{normal form} if there are r degree 1 (i.e., completely multiplicative) factors and s inverses of degree 1 factors, and if no two factors in the product are mutually inverse to one another. By commutativity we can always write the positive (degree 1) factors first.
\end{definition}
 
\begin{theorem}\label{mytheorem10}
Let $\alpha\in$ \textbf{MF}
 \vspace{0.5cm}

(1) $\alpha$ is type 1 = $(fin,fin)$ if and only if it is the identity.
 \vspace{0.5cm}

(2) if $\alpha$ is positive  it is type 2, $(\infty ,fin)$; if $\alpha$ is negative, it is type 3 .
 \vspace{0.5cm}
\end{theorem}  
 \textbf{Proof} 

 (1) Suppose that $\alpha$ and $\alpha^{-1}$ have a finite range, i.e., both  $\alpha(p^n)$ and $\alpha^{-1}(p^n)$ have non-zero values for infinitely many values of $n$.  Let  $F_n$
represent $\alpha$ and $F'_n$ represent  $\alpha^{-1}$.  Then let $F'_j = t_j,  j=1,...,k, F'_n =0,$ otherwise.  We use the linear recursion properties of the GFP's,  and suppose that $F_n\neq0, n>0 $,  but that $F_{n+1} = F_{n+2}=...=F_{n+k} =0$ (a condition that is sufficient to guarantee that all of the $F_j = 0, j > n$).  Then  we have the following relation $$0 = F_{n+k} = t_1F_{n+k-1} +...+t_{k-1}F_{n+1} + t_kF_n.$$ But all of the terms on the right hand side are equal to  $0$ except $t_kF_n$,  which is assumed to be non-zero. This yields a contradiction.  The case $F_n=0, n>0 $ is just the case when the arithmetic function is the convolution identity, $F_0=1, F_n=0, n>0 $. Thus we have proved that there are three mutually exclusive types,  and that type 3 contains only the identity. 

(2)  claims that  functions with valence $< r,0>$ have infinite ranges and finite domains;  while functions with valence  $<0,s>$  have finite ranges and infinite domains.    But   these two propositions follow from the fact that we have defined functions to be positive if  they are a product of CM functions, and negative if they are the inverses of a product of CM functions, from Theorem \ref{mytheorem5}, and from  the fact that the core polynomials determine the number of parameters of a function.      $\square$                                                                                                                                                                         

\begin{remark}\label{myremark2}
We shall show that $\phi $ is of type 4. 
To show that $\phi$  is of type 4, we must show that both $\phi$ and  $\phi^{-1}$  have infinite non-zero range. But,  for any prime  $p$ and all $n$, $\phi{(p^n)} = p^n-p^{n-1} \neq 0$. $\phi^{-1}$ is also infinite:  its  parameters are $t_j = -(p-1)$ for all $n$ .  
So let us suppose that $\overline{F}_j = 1-p$ for  $0<j<k$.  Then,  by the recursive property of GFP-sequence and induction ,  we have that
$$\overline{F}_k = F_1\overline{F}_{k-1}+ ...+F_{k-1}\overline{F}_1 + F_k $$
 $$= (p-1)(F_1+F_2+...+F_{k-1} )-F_k  $$
 $$= (p-1)(\sum_{j=1}^{k-1} F_j) - F_k $$
 $$=(p-1)(\sum_{j=1}^{k-1}(p^j-p^{j-1})) -p^k+p^{k-1}$$
 $$=(p-1)(p^{k-1}-1)-p^k+p^{k-1}$$
  $$=1-p$$
  Thus we have shown that $t_j =1-p\neq 0$ for all j,  and so,  $\phi\in$ type 4. $\square$

  We can also prove this by using the fact that  $\phi = \zeta_1 \ast \mu$.  Represent  $\zeta_1$ by $F'$ and $\mu $ by $F''$, then by Corollary \ref{mycorollary9}, $ t'_1 = p, t'_j = 0, j>1;  s'_n = -p^n$; $s''_n = -p,$ when $n=1,$ and $0$ otherwise, and $t''_n = 1$ for all values of $n$.  We can symbolize the product by $(\infty,finite)\ast (finite, \infty) $ .
\end{remark}

This property of  $\phi$ turns out to be true for all totient functions, i.e., functions of valence $<1,1>$.

\begin{proposition}\label{myproposition7}
$\alpha = \beta \ast \gamma^{-1}$  where $\beta, \gamma$ are degree 1 functions,   $\alpha$  is not the  identity function then $\alpha$ is type $(\infty , \infty)$ and degree $\infty$.
\end{proposition}

\textbf{Proof}.
We represent $\alpha, \beta,  and \; \gamma$ by, respectively, $F,F', and \: F''$, with $t_j, t'_j,t''_j; s_j,s'_j,s''_j$ as in Corollary 9.  We observe that $F'_n =t'^n_1, and \: F''_n =t''^n_1$, and that neither of these values is 0.  We also have that $F_n = t'^n_1 - t'^{n-1}_1t''_1$, since $t''_j = 0, j>1$,  and $F_n = 0$,
implies $t'_1= t''_1$, contradicting the hypothesis;  hence, the range of  $\alpha$ is infinite.  If we apply similar reasoning to $F_n^{-1} $, we find that it also has infinite range,  thus $\alpha$ is of type 4 and has infinite degree. $\square$

\vspace{0.5cm}
In a paper in 2005 $\cite{LP}$,  Laohakosol and Pabhapote discussed Busche-Ramanujan identities and the Kesava Menon norm. We shall discuss the Kesava Menon norm in the next section.  Now we wish to look at their theorem extending Busche-Ramanujan identities to multiplicative functions of mixed type, Corollary 2.4 in \cite{LP},  which we reproduce here:

Let  $\alpha \in MF,$ Then the following hold. 

(i)  $\alpha \in \mathcal{C}(1,1) \Longleftrightarrow  $ for each prime $p$ and each  $n \in \mathbb{N},$ there exists a complex number $T(p)$ such that $$ \alpha (p^n) = T(p)^{n-1} \alpha (p).$$

(ii) $\alpha \in \mathcal{C}(2,0) \Longleftrightarrow $ for each prime $p$ and each $ n (\geqslant 2) \in \mathbb{N},$ $$\alpha (p^{n+1}) = \alpha(p) \alpha(p^n) + \alpha (p^{n-1})[\alpha (p^2)-\alpha(p)^2].$$

(iii) $\alpha \mathcal \in {C}(1,s) \Longleftrightarrow$ for each prime $p$ and each $\alpha \in \mathbb{N},$ there exist complex numbers $B_1(p),...,B_s(p)$ such that for all $\alpha \geqslant s,$ $$\alpha (p^n) = \sum_{j=0}^s \gamma(p)^{n-j}H_j$$
 where $$H_j = (-1)^j \sum_{1\leqslant i_1 <i_2 < ...i_j \leqslant s } B_{i_1}(p)...B_{i_j}(p), \; H_0 = 1.$$
 $\gamma \in CM$.
\vspace{0.25cm}

We have made only minor changes in the notation in keeping with the notation established in this paper. Some of the functions mentioned here were explained in the context of a theorem which preceeded this corollary in the cited article.  Features of this theorem will be mentioned below, but for the time being it is not necessary to clarify these statements further.  Rather we wish to apply the methods of our paper to show how these results can be clarified and simplified.

First we note that (ii) is just McCarthy's theorem  \cite{PJM} discussed in section 5, where we showed that the mysterious bracket is just the parameter $t_2$ of the positive degree 2 multiplicative function in question. We look now at parts (i) and (iii) of the Laohakosol-Pabhapote corollary.  We shall, as is consistent with our practice in this paper, drop reference to the prime $p$ since our theory is instrinsically local.  So we now look at part (i).

Let $\alpha =  \beta \ast \gamma^{-1} $ where $\alpha, \beta, \gamma \in CM$.  (Recall the discussion of the Euler totient function $ \varphi$ above.) Let  $F(\textbf{t}) \leftrightarrow \alpha, F'(\textbf{t'}) \leftrightarrow\beta, F"(\textbf{t''}) \leftrightarrow \gamma$. Then $$t'_1 \neq 0, t'_j = 0, j >1; t''_1\neq 0, t''_j = 0, j >1 $$
Making use of Proposition \ref{myproposition6}, we have that $$ \overline{s}" = \overline{F}"_n = -t"_1 \, if \,n=1, = 0 \,if\, n>1,$$

From this we easily deduce that:
$$F_1 = t'_1-t''_1 = t_1,$$ $$ ...,$$ $$ F_n = t'^{n-!}_1 (t'_1-t''_1) = t'^{n-1}_1(t'_1-t''_1) =  t'^{n-1}_1t_1 =  t'^{n-1}F_1$$
That is, 
\begin{proposition}\label{myproposition8}

Suppose  $\alpha= \beta \ast \gamma^{-1}$, $\beta \; and \; \gamma$ degree 1 functions, $ \beta \neq  \gamma$ and $ \beta \neq \delta \neq \gamma$, and suppose that  $\alpha, \beta,  and \; \gamma$ are represented by, respectively, $F,F', and \: F''$, with  parameters $t_j, t'_j,t''_j; s_j,s'_j,s''_j$.  Then, 
$$F_n =  t'^{n-1}_1F_1, \; \forall n\in\mathbb{N}$$
where $ \overline{F}'' = F^{-1}.$ $\square$
\end{proposition}

Thus, the mysterious $T(p)$ in the original corollary is just one of the parameters determining the representing GFP-sequence just as in the case of McCarthy's theorem, this time for $\beta$, it is just $t'_1$. 

Next,  we look at part (iii) of the corollary.  We first consider a product of degree 3, as in Remark 3, i.e., a product $\alpha = \beta  \ast  \gamma^{-2}$  where $\beta$ is a positive function of degree 1, and  $\gamma$ is a positive function of degree 2, and use the formalism of the previous example for representing  $\alpha, \beta$ and $\gamma$ by GFP sequences.  Thus we have $$F(t_1, t_2,t_3) = F'(t'_1)\ast \overline{F}''(t''_1,t''_2).$$  (Here we are thinking of $ \overline{F}''$ not  as an inverse, but as a function in its own right ).  With the help of Proposition \ref{myproposition6} and Corollary \ref{mycorollary9}, we easily find that 
$$F_1 = t'_1+t"_1 = t_1$$ 
and 
$$F_2 = F'_2 + F'_1\overline{F}''_1+\overline{F}''_2 =  t'^2_1+t'_1t''_1+t''_2$$ 
$$= t'_1(t'_1+t''_1) = t'_1F_1-s''_2.$$ 

In the same way, $$F_3 = t'^2_1F_1-t'_1s''_2-s''_3.$$

Induction gives us
\begin{proposition}\label{myproposition9} 
$\alpha \in <1,s> $  $\alpha \leftrightarrow F$
$$F_n = t'^{n-1}F_1- t'^{n-2}_1s''_2-t'^{n-3}_1s''_3. \quad\square$$
\end{proposition}

The proposition is nothing but a thinly disguised version of the definition of convolution product together with the particular assumptions about the parameters of the factors together with Proposition \ref{myproposition6}.  But it says all that the Laohakosol-Pabhapote result says and at the same time explicitly identifies the mysterious functions.  Compare the previous remarks concerning the McCarthy theorem. The general case is now clear.

\begin{theorem}\label{mytheorem11}
Let  $\alpha =\beta \ast \gamma$,  where  $\beta$ is a positive degree 1 function, and  $\gamma$ is the convolution inverse of a degree $(k-1)$,  positive function, and suppose that 
$\alpha$ is non-trivial and in normal form.  Representing these functions by GFP functions as above $\alpha \leftrightarrow F(t_1,...,t_j,...),\beta \leftrightarrow F'(t'_1), \gamma \leftrightarrow F''(t''_1,...,t''_j,...)$, 

$$F_{k,n} = -\sum_{j=0}^n t'^{(n-j)}_1s''_j =  \sum_{j=0}^n t'^{(n-j)}_1 F''_j$$
\end{theorem} 

\textbf{Proof}
$\beta \leftrightarrow F'$ is of type 2, and $\gamma \leftrightarrow F''$ is of type 3, thus as a result of the assumptions on the factors and the freeness of the product, $F$ is infinitely generated, and so is its inverse;  thus $F$ is type 4;  we can symbolize this by  $(\infty,fin)\ast (fin,\infty)$.  
By Corollary \ref{mycorollary9},  we have that $$F_n = -s'_n + \sum_{j=1}^{n-1} s'_{n-j}s''_j -s''_n,$$
and since $\beta$ is of degree 1,  $$-s'_j = F'_j = t'^j_1$$
Therefore,
$$F_n = -\sum_{j=0}^n t'^{(n-j)}_1s''_j ,$$
That is,  $$F_n = \sum_{j=0}^n t'^{(n-j)}_1F''_j $$ $\square$

Since $\alpha$ in Theorem \ref{mytheorem11} has valence $<1,r>$, Theorem \ref{mytheorem11} is a generalization of  Proposition \ref{myproposition8}.   

The following corollary is a direct consequence of the remarks in the proof of the previous theorem concerning the generators and the types of terms in the product:
\begin{corollary}\label{mycorollary12}
A non-trivial product in normal form of valence $<r,s>$, that is, a \textit{ mixed} product, is type 4.
\end{corollary}

\vspace{0.05cm}
In \cite{LP} the notion of \textit{s-excessive} was introduced in connection with a generalization of the Busche-Ramanujan identities. Since the theory that we are dealing with is local,  we shall refer the reader to Definition 3.1 of the paper just cited for the general definition and discussion of that concept and give here the local definition, which is suitable for this paper.

\begin{definition}\label{mydefinition7}
Two prime powers of the same prime $p$, say $p^r$ and $p^s$  with $s \leqslant r$, are said to be \textit{e-excessive} if $r-s=e$.  
\end{definition}

In section 5,  we suggested that Busche-Ramanujan identities can be regarded locally as expressing  $F_n = F_{r+s}$ in terms of $F_r$ and $F_s$. The following theorem,Theorem \ref{mytheorem13},  is a local generalization of B-R identities to functions of arbitrary degree,  and generalizes to Theorem 3.2, \cite{LP} when the functions here are restricted to the functions in that theorem and the result of Theorem \ref{mytheorem13} is globalized. 

\begin{theorem}\label{mytheorem13}
Let $\alpha$ be a multiplicative arithmetic function of degree $k$, and let $r$ and $s$ be two integers with $r \leqslant s$,  and abbreviating $F_{k,n} $ as $F_n$, then $$F_n = F_{r+s} = \sum_{j=0}^{e+1} (-1)^j S_{(r,1^{j})} F_{s-j}.$$
where $S_{(r,1^{j})}$ is an isobaric reflect of the Schur-hook function whose Young diagram has an \textit{arm} of length $r$ and a \textit{leg} length of $j$. 
 \end{theorem}
 Observe that these Schur-hook functions for a given $r$ consist exactly of the elements of the row $r-th$-row, in order from right to left, of the companion matrix, and that $ S_{(r,1^{0})}= F_r$;  so that this formula satisfies our interpretation of a generalization of the Busche-Ramanujan identity,  and includes Theorem 3.2 in \cite{LP}. Also observe that such a row is a vector representing the $r$-th power of a root of the core polynomial (see section 2). 
\vspace{0.05cm}

 \textbf{Proof} (of Theorem \ref{mytheorem13})
First note that the companion matrix is \textit{stable} in the sense that adding a $k+1-$st column on the left of an $\infty \times k-$companion matrix changes nothing in the original matrix.  Thus we might as well assume that we have infinitely many parameters for $F$. For a finite $k$ we need only let $t_j = 0$ from a certain point onward.  We shall also need to recall the fact that the columns of the companion matrix are linear recursions with respect to the parameters $t_j$. When $r = 0$ the theorem is just a statement of the fact that  $F$ is a linear recursion with parameters $t_j$.  We proceed by induction on $r$ with $e=r-s$. 
$$F_n = F_{r+s} =  t_1F_{(r-1)+s}+t_2F_{(r-2)+s}+...+t_kF_{(r-k)+s}$$

$$= t_1\sum_{j=0}^{e+2}S_{(r-1,1^j)}F_{s-j}+t_2\sum_{j=0}^{e+3}S_{(r-2,1^j)}F_{s-j}+...+\sum_{j=0}^{e+k+1}S_{(r-k,1^j)}F_{s-j}$$
$$=(\sum_{j=0}^{e+2}t_1S_{(r-1,1^j)} +\sum_{j=0}^{e+3}t_2S_{(r-2,1^j)}+...+\sum_{j=0}^{e+k+1}t_kS_{(r-k,1^j)}) F_{s-j}$$
$$=\sum_{j=0}^{e+1}S_{(r,1^j)}F_{s-j} \quad\quad \square$$

Corollary 2.6 \cite{LP} simply states that for a positive function of degree $r, \; t_j = 0, j>r, t_r \neq 0$. And  
Corollary 2.7 \cite{LP} is a statement of the basic fact that the representing GFP-sequence is a k-order linear recursion whenever $\alpha$ is a multiplicative function of degree $k$, \cite{MW} or, indeed, \cite{TM2}, \cite{MT}, and \cite{MT2}.

To conclude this section we point out that some  \textit{classical} theorems in multiplicative number theorem reduce to rather prosaic statements, if not trivialities, in the formalism of this paper, that is in the context of the theory of symmetric functions. An example of this is the binomial identity, \cite{ReS}, \cite{RS2},\cite{PH}.   We do not wish to labour this point,  but we would like to use this instance to show how this is the case.

The binomial identity may be stated as follows:

\begin{proposition}\label{myproposition10}
Suppose that  $\alpha = \gamma_1 \ast \gamma_2$ where $ \gamma_1$ and  $\gamma_2$ are completely multiplicative functions, then $\alpha$ satisfies the binomial identity
$$\alpha(p^n) = \sum_{j=0}^{\lceil\frac{n}{2}\rceil}(-1)^j\left(\begin{array}{cc}n-j &  \\j & \end{array}\right)\alpha(p)^{n-2j}(g_1(p)g_2(p))^j, \quad\quad \square$$
\end{proposition}

\noindent which in turn is just a special case of the general formula for Weighted Isobaric Polynomials, \cite{MT}.
If we use our GFP representation,  then this result is just a case of the general formula for GFP's as in section 1 of this paper; namely, 

$$F_{k,n} = \sum_{\alpha \vdash n} 
  \left( \begin{array} {c}   | \alpha | \\  \alpha_1,...,\alpha_k  \\   \end{array}
\right)  t_1^{\alpha_1}... t_k^{\alpha_k}  . $$

Here is the formulation in our terms:
\begin{proposition}\label{myproposition11}
Let $\alpha$ be a positive MF of degree 2 represented by $F$, and let $F'$ and $F''$ represent the convolution factors, then $$F_{2,n} = F_n = \sum_{j=0}^{\lceil\frac{n}{2}\rceil}(-1)^j \left(\begin{array}{cc}n-j &  \\j & \end{array}\right) F_{n-2j}(-t'_1t''_1)^j.   \quad\quad \square$$
\end{proposition}

Thus, we have $$F_n = \sum_{j=0}^{\lceil\frac{n}{2}\rceil}(-1)^j \left(\begin{array}{cc}n-j &  \\j & \end{array}\right) F_{n-2j}(-t_2)^j.$$ 

Letting $n=5$,  and noting that the partitions of $5$ with $k=2$ are just $(1^5),  (1^3,2), (1,2^2)$,
we have, $$F_{2,5}= \sum_{j=0}^{3}(-1)^j \left(\begin{array}{cc}5-j &  \\j & \end{array}\right) F_{5-2j}(-t_2)^j= t_1^5+4t_1^3t_2+3t_1t_2^2 .$$ 
\textbf{8. NORMS}
\vspace{0.5cm}

There are two norms that are of importance in the theory of arithmetic functions.  $\cite{R}, \cite{HS}$.  The first of these gives as the norm of an AF  $f$,  the least integer $n$ such that  $f(n)\neq0$.  In this norm all non-trivial MF have norm $1$.  So this is not a useful norm for our purposes.  

The second norm $\mathcal{N }$ is defined only on multiplicative functions and is given by $$\mathcal{N}(\alpha)(n) =  \sum_{d|n^2} \alpha(n^2) \lambda(d) \alpha(d)$$
where  $\lambda$ is the Liouville function.  $\mathcal{N }$ is  a multiplicative function.  If  $\alpha$ has degree 1 or 2,  then $deg\mathcal{N}(\alpha)=1,2$, respectively. ($\cite{PJM}, p.50$)
What does the norm look like in isobaric notation?
\vspace{0.5cm}

Let $F_n(\textbf{t}) \leftrightarrow  \alpha, \alpha\in MF$ for suitable choice of the vector  $\textbf{t}$, and let  $\mathcal{N}(F_n)= \mathcal{N}(\alpha)(p^n)$.
Then $$\mathcal{N}(F_n) =  \sum_{j=0}^{2n} (-1)^j F_{2n-j} F_j$$
or, equivalently
$$\mathcal{N}(F_n) = (2 \sum_{j=0}^{n-1} (-1)^j F_{2n-j} F_j) + (-1)^nF_n^2$$

By Corollary \ref{mycorollary9}, this can be written in the following way.

\vspace{0.05cm}
\begin{proposition}\label{myproposition12}
$$\mathcal{N}(F_n) = -s_{2n} + 2\sum_{j=1}^{n-1} (-1)^j s_{2n-j}s_j+ (-1)^ns_n^2$$
$$\mathcal{N}(F_n) = -t_{2n} + 2\sum_{j=1}^{n-1} (-1)^j t_{2n-j}t_j+ (-1)^nt_n^2. \;\; \square $$
\end{proposition}

It is well known that this norm is multiplicative, that is, that: If  $\alpha$   and $\beta \in$ MF,  then  
$$\mathcal{N}(\alpha \ast \beta)=\mathcal{N}(\alpha) \ast \mathcal{N}(\beta)  (\cite{PJM}, p.50).$$

It is useful to give a proof of this fact here using Theorem \ref{mytheorem1} .

\begin{theorem}\label{mytheorem14}
Let $F = F'\ast F'',$ then $$\mathcal{N}(F_n) = \mathcal{N}(F'_n) \ast \mathcal{N}(F''_n);$$
that is, if $\mathcal{N}_n = \mathcal{N}(F_n),$ then $$\mathcal{N}_n = \mathcal{N}'_n \ast \mathcal{N}''_n$$
\end{theorem}

\begin{lemma}\label{mylemma2}
Let the set of parameters for the MF $\mathcal{N}$, where $\mathcal{N}$ is the norm function,  be denoted by $n_r$.  Then $$n_r =  n'_r - \sum_{j=1}^{r-1}n'_jn''_{k-j} +n''_k.$$
\end{lemma}

\noindent \textbf{Proof} of Lemma.
\vspace{0.05cm}

This is just an application of Corollary \ref{mycorollary9}. $\square$ 

\noindent \textbf{Proof} of Theorem. 

Using the recursion property of the norm function and the lemma we have

$$\mathcal{N}_r = \sum_{j=1}^r n_j \mathcal{N}_{r-j} $$
and

$$\sum_{j=1}^r (n'_j - \sum_{s=1}^{j-1}n'_{j-s}n''_s +n''_j)\mathcal{N}_{r-j}.$$

Induction gives us that

$$= \sum_{j=1}^r (n'_j - \sum_{s=1}^{j-1}n'_{j-s}n''_s +n''_j)(\mathcal{N'}_{r-j} \ast \mathcal{N''}_{r-j}) $$ 

which in turn is

$$= \sum_{j=1}^r (n'_j - \sum_{s=1}^{j-1}n'_{j-s}n''_s +n''_j)(\sum_{i=0}^j(-1)^i\mathcal{N'}_{r-j-i} \mathcal{N''}_{i}) $$

$$= \sum_{j=1}^r n'_j\mathcal{N}_{r-j}-( \sum_{s=1}^{j-1}n'_{j-s}n''_s)\mathcal{N}_{r-j}+n''_j\mathcal{N}_{r-j}.$$

This can be rewritten as
$$ [(\sum_{j=1}^r n'_j \mathcal{N}'_{r-j}  + \sum_{j=1}
^r n''_j \mathcal{N}''_{r-j} ) + ... +(\sum_{j=1}^r n'_{j-1} \mathcal{N}'_{s-(j-1)} )\mathcal{N}''_{r-s} +( \sum_{j=1}^r n''_{j-1} \mathcal{N}''_{s-(j-1)})\mathcal{N}'_{r-s} ] $$
 $$ +[ \sum_{s=1}^{r-1} (n'_s(\mathcal{N}'_{r-s}-\sum_{j=1}^{r-s}n'_j \mathcal{N}'_{r-s-j}) +n''_s(\mathcal{N}''_{r-s} - \sum_{j=1}^{r-s}n''_j \mathcal{N}''_{r-s-j}))]$$
writing this as $[A] +[B]$, and noting that  $[B] = 0$ and $A = \sum_{j=0}^r \mathcal{N}'_{r-j}\mathcal{N}''_j $. we have that the above expression 

$$ =  \sum_{j=0}^r \mathcal{N}'_{r-j}\mathcal{N}''_j  $$
$$=  \mathcal{N}'_r \ast \mathcal{N}''_r , \quad\quad \square$$

\begin{corollary}\label{mycorollary15}
Let \{F\}, \{F'\},\{F''\} represent  multiplicative functions locally and suppose that $F_1 = F'_1 + F''_1$,  then
$F = F' \ast F''$ $\square$
\end{corollary}

\begin{lemma}\label{mylemma3}
If $deg(\alpha) = 1$, then  $deg\mathcal{N}(\alpha) = 1$.  ($\cite{PJM}, p.50)$
\end{lemma}

\begin{theorem}\label{mytheorem16}
$$deg(\mathcal{N}(\alpha)) = deg(\alpha)$$
\end{theorem}
\textbf{Proof}   If $deg(\alpha) =k$, then $\alpha$ is the convolution product of  $k$ degree 1 multiplicative functions.  The theorem then follows from the two lemmas. $\square$

\vspace{0.5cm}
		
SOME EXAMPLES:

1.  Let $\alpha$ be the arithmetic function whose values for $\alpha (p^n) =f_n, \forall $ primes $p$, where $f_n$
is the $n-th$ Fibonacci number.  $f_0=f_1 = 1, f_{n+1} = f_n + f_{n-1}$, and represent $\alpha$ at each prime by $F$.  It is easy to calculate that $F$ has degree two. In fact, that $F_n(t_1,t_2)= F_n(1,1)$. According to Theorem \ref{mytheorem16},  $deg\mathcal{N}(\alpha) = 2$ also.  We can see this directly.  $\mathcal{N}(F_1) = \mathcal{N}_1 = f_3 = 3 $ and $\mathcal{N}_2 = f_5 = 8$.  From the first of these two facts, we have that $n_1 = f_3 = 3$ and from the second, that $\mathcal{N}_2 = f_5 = 8$, and since $\mathcal{N}_2 = n_1^2 + n_2$,  we have that $n_2 = -1$.  The same technique shows that $\mathcal{N}_3 = 0$. An easy induction yields $n_j = 0, j>2$.  We can then use the recursion property to show that $\mathcal{N}_n = f_{2n+1}$.

2.  For the multiplicative function  $\tau$ (number-of-divisors function), using Theorem \ref{mytheorem16} and the methods above,  we find that $deg\mathcal{N}(\tau) = 2 = deg(\tau)$, where $t_1 = 1, t_2 = -1$ and $\mathcal{N}_n = n+1 = \tau(p^n)$ $n_1 = 2, n_1 = -1$.		

\vspace{0.5cm}
9. \textbf{WIP-MODULE REVISITED}

\vspace{0.5cm}
In Section 1, the $\mathbb{Z-}$module generated by the Schur-hook polynomials was introduced.  In Section 2, the infinite companion matrix $A^\infty$ was introduced.  Among the properties of this matrix discussed there, it was mentioned that the columns were sequences of Schur-hook functions,  the right-most being the sequence of generalized Fibonacci polynomials. If the columns are numbered from right to left as $0,1,...,k-1$, the  $j-th$ column in this numbering is the sequence $\{(-1)^j S_{n,1^j}\}_n$, 
where $n$ varies over all of the integers, positive and negative and zero.   In particular this extends the definition of any MF to negatively indexed functions.  Moreover,  each of the columns represents a MF locally; in fact, any linear combination of sequences does so.  Thus a given core polynomial (of finite or infinite degree) picks out a specific MF, the one determined by the GFP-sequence; but at the same time, the companion matrix brings to life a basis for an entire module of MF's associated with the same core.   This phenomenon of associated MF's should be interesting to investigate. 

In an earlier work which appeared in \cite{MT} (also see \cite{C} and \cite{CG} ), a construction was given which provided a $q$-th root, $q\in\mathbb{Q}$,  for every element of the WIP-module, thus giving an effective construction of the divisible closure (injective hull) of the WIP-module.  So, in particular, providing a construction for the $q-$roots of any multiplicative function.  We shall repeat this formula here.  In \cite{CG} Carroll and Gioia gave a numerical description of the rational roots of the elements of the group generated by the completely multiplicative functions $\mathcal{L}$. In \cite{MT}, a collection of isobaric polynomials with rational coefficients was found which play the same role for the rational roots of the multiplicative arithmetic functions as the GFP's play for the group $\mathcal{L}$ itself, thus embedding this group into a divisible group.  In fact, we do more:  we provide each member of the WIP-module with a unique rational root induced by convolution. Taking the injective limit then gives an embedding of $\mathcal{L}$ in its divisible closure analogous to the construction of the rationals from the integers.  The main theorem  (Theorem 5.1 in \cite{MT}) is the following.

\vspace{0.25cm}
\begin{theorem}\label{mytheorem17}
Let $H_{k,n}(t,q)$ denote the $q-th$ convolution root of  $F_{k,n}(t)$,  where
 $F_{k,n}(t) \in GFP$, then
$$H_{k,n}(t,1) \sum_ {\alpha \vdash n} \frac{1}{|\alpha |!} B_{(|\alpha-1|)}^q  \left(\begin{array}{ccc} & |\alpha| &  \\\alpha_1 & ... & \alpha_k\end{array}\right)t^{\alpha}$$
is the $q-th$-convolution root of the weighted isobaric polynomial  $P_{k,n,\omega}$. 
In particular,  $H_{k,n}(\textbf{t},1)  = F_{k,n}(\textbf{t}),$  (Corollary 5.2,  \cite{MT}). $\quad\quad \square$  
\end{theorem} 
\vspace{0.05cm}

For a discussion and proofs see section 5 of  \cite{MT}. The consequences for AF are that every multiplicative arithmetic function in the Dirichlet ring is now equipped with a unique q-th convolution  root for every rational number q. Moreover,  we have an embedding of $\mathcal{L}$ into its divisible closure.
\vspace{0.05cm} 

In another paper \cite{MW} in which the authors investigated periodicity with respect to linear recursions, we were able to apply these ideas to algebraic number fields. The consequences of the results in that paper for arithmetic functions are:
\begin{theorem}\label{mytheorem18}
An arithmetic function is periodic if and only if the roots of its core are all roots of unity.
In particular, the period of a cyclotomic core is just the order of the associated cyclotomic group.
\end{theorem}

\begin{theorem}\label{mytheorem19}
For every prime number $p$, every MF induced by a core of degree $k$ is periodic modulo-p with a period dividing $p^k-1$.
\end{theorem}\vspace{0.05cm}

\begin{remark}\label{myremark3}
 A computer program for finding this period is given at the end of \cite{MW}.
\end{remark}

In \cite{MW}, the theory of periodicity was used to study algebraic number fields.  Thus we have a triad of related ideas,  all linked together through the theory of symmetric functions:  namely,  linear recursions, number fields, and multiplicative arithmetic functions..

In a later paper we shall discuss the relation of the subgroup $\mathcal{L}$ to the group $\mathcal{M}$, as well as  Ramanujan sums and the Ramanujan tau-function in terms of the Isobaric representation.

  Key Words: multiplicative arithmetic functions, symmetric polynomials,  linear recursions.
  O5EO5,11B37,11B39,11N99
  
  \vspace{0.5cm}
  
 \noindent Trueman MacHenry \\
  York University, Toronto, Canada \\
  machenry@mathstat.yorku.ca \\
  Kieh Wong \\
  Centennial College, Toronto, Canada \\
  kkwong@centennialcollege.ca
 
 \end{document}